\documentclass{article}

\usepackage{amsfonts}
\usepackage{xypic}

\title{Virtually Free pro-$p$ groups whose Torsion Elements have finite Centralizer}%

\author{W. Herfort\thanks{The first author would like to thank for the hospitality
  at the UNC during February 2003}, P.A. Zalesski\thanks{The second author expresses his thanks for support through CNPq}
}

\newif\iftest % Testversion / Finalversion
\testfalse %%%%%%%%%%%%%%%%%%%%%%%%%%%%%%%%%%%%%%%%%%%%%%%%TESTTEST
\newif\ifteston
\iftest\testontrue\else\testonfalse\fi
\newif\ifnotes % notes enabled
\notestrue
\newif\iffinal 
\finalfalse
\iffinal\global\notesfalse\global\testfalse\fi  % 
\newcounter{notes}
0

\newcommand{\ignoriere}[1]{}

\def\plabel#1{\label{#1}\iftest\fbox{#1 }\fi}
\def\pref#1{\ref{#1}\iftest{ \fbox{#1 }}\fi}
\def\pcite#1{\cite{#1}\iftest{\fbox{#1 }}\fi}  
\newcommand{\prefeq}[1]{Eq.(\pref{e-#1})}

\newcommand{\preflemma}[1]{Lemma \pref{l-#1}}
\newcommand{\prefdef}[1]{Definition \pref{d-#1}}
\newcommand{\prefprop}[1]{Proposition \pref{p-#1}}
\newcommand{\preftheorem}[1]{Theorem \pref{t-#1}}

\newcommand{\prefrem}[1]{Remark \pref{rem-#1}}                             

\newtheorem{theorem}{Theorem}[section]
\newtheorem{lemma}[theorem]{Lemma}
\newtheorem{definition}[theorem]{Definition}

\newtheorem{proposition}[theorem]{Proposition}

\newtheorem{corollary}[theorem]{Corollary}
\newtheorem{remark}[theorem]{Remark}

\newcounter{claims}

\newcommand{\bd}[1]{\begin{definition}\plabel{d-#1}\rm}
\newcommand{\ed}{\end{definition}}
\newcommand{\bt}[1]{\begin{theorem}\plabel{t-#1}\setcounter{claims}0}
\newcommand{\et}{\end{theorem}}
\newcommand{\bl}[1]{\begin{lemma}\plabel{l-#1}\setcounter{claims}0} 
\newcommand{\el}{\end{lemma}}               
\newcommand{\bc}[1]{\begin{corollary}\plabel{c-#1}}
\newcommand{\ec}{\end{corollary}}
\newtheorem{notation}{Notation}
\newcommand{\bn}[1]{\begin{notation}\plabel{n-#1}\rm}
\newcommand{\en}{\end{notation}}
\newcommand{\brem}[1]{\begin{remark}\plabel{rem-#1}\rm}
\newcommand{\erem}{\end{remark}} 
\newcommand{\bp}[1]{\begin{proposition}\plabel{p-#1}\setcounter{claims}0}
\newcommand{\ep}[1]{\end{proposition}}
\newcommand{\be}[1]{\begin{equation}\plabel{e-#1}}
\newcommand{\ee}{\end{equation}}
\newlength{\help}
\newcounter{subclaims}
\newenvironment{env}{}{}

\newcommand{\bcl}{\medskip\addtocounter{claims}1\setcounter{subclaims}0{
\noindent{\it Claim \arabic{claims}: } }
                  \begin{env}\it}

\newcommand{\ecl}{\end{env}\rm\medskip}
\newcommand{\bscl}
  {\addtocounter{subclaims}1{Subclaim \arabic{subclaims}}\begin{env}\it}
\newcommand{\escl}{\end{env}\rm\smallskip}

\newcommand{\claimproof}{{}}

\newcounter{inrmlist}
\newcounter{inalphlist}
\newenvironment{rmenumerate}{
                       \setcounter{inrmlist}{0}
                       \begin{list}
                         {\makebox[2em]{(\roman{inrmlist})}}{
                         \addtolength{\leftmargin}{1em}
                         \addtolength{\itemsep}{1mm}
                         \addtolength{\itemindent}{-2em}\usecounter{inrmlist}}
                       }%
                       {\end{list}}

                       {\end{list}}

\newcommand{\Z}{{\bf Z}}

\newcommand{\of}[1]{{(#1)}}
\newcommand{\tF}{\tilde F}
\newcommand{\tG}{\tilde G}

\newcommand{\tv}{\tilde v}
\newcommand{\tw}{\tilde w}

\newcommand{\cdc}{contradiction}
\newcommand{\cex}{counter-example}
\newcommand{\epi}{epimorphism}
\newcommand{\fg}{finitely generated}
  
\newcommand{\fpp}{free \pp}

\newcommand{\fppgrp}{free \ppgrp}

\newcommand{\pp}{pro-$p$}

\newcommand{\ppgrp}{\pp\ group}

\newcommand{\vfpp}{virtually \fpp}
\newcommand{\vfppgrp}{\vfpp\ group}

\newcommand{\path}[2]{[#1,#2]}                %  path beginning at 1 ending at 2  

\newcommand{\inv}{^{-1}}                        % inverse of an element
\newcommand{\indx}[2]{\vert #1 : #2\vert}          %  index of a subgroup in a group

\newcommand{\rank}{\text{rank}}
 
\newcommand{\tor}[1]{\text{Tor}(#1)}             % Torsion 
\newcommand{\gp}[1]{\langle #1 \rangle}             % Group generated by #1
\newcommand{\torgp}[1]{\gp{\tor{#1}}}              % Group generated by torsion
\newcommand{\torfactor}[1]{#1/\torgp{#1}}           % factorgp G/ tor G 

\newcommand{\ol}[1]{\overline{#1}}
\newcommand{\ugp}{\{1\}}                          % Trivial group
\def\xtimes{\times\hskip -0.2em\vrule height 4.6pt depth -0.3pt\hskip0.2em}
\newcommand{\rtimes}{\mbox{$\xtimes$}}
\newcommand{\multip}{g_{\tw}\inv}

\newcount\names
\def\aname#1#2{\ifodd\names{#1.\,#2}\else{#2 #1.}\fi}

\newcommand{\efrat}{\aname{I}{Efrat}}

\newcommand{\haran}{\aname{D}{Haran}}
\newcommand{\herfort}{\aname{W}{Herfort}}

\newcommand{\melnikov}{\aname{O.V}{Mel'nikov}}

\newcommand{\ribes}{\aname{L}{Ribes}}

\newcommand{\weiss}{\aname{A}{Weiss}}
\newcommand{\scheiderer}{\aname{C}{Scheiderer}}

\newcommand{\segal}{\aname{D}{Segal}}

\newcommand{\shalev}{\aname{A}{Shalev}}

\newcommand{\zalesskii}{\aname{P.A}{Zalesskii}}

\def\HNN{{\it HNN}}
\def\HNNgrp{HNN-group}
\def\HNNext{HNN-extension}

\def\KST{Kurosh subgroup theorem}

\def\minass{minimality assumption}

\newcommand{\homo}{homomorphism}

\def\KST{Kurosh subgroup theorem}

\begin{document}
\maketitle

\begin{abstract}
A \fg\ \vfppgrp\ with finite centralizers of
its torsion elements is the \fpp\ product of finite $p$-groups and
a \fpp\ factor.
\end{abstract}

\testfalse

\section{Introduction}

The objective of this paper is to give a complete description of a
finitely generated  virtually free pro-$p$ group whose torsion
elements have finite centralizers. Our main result is the
following

\bt{main} Let $G$ be a \fg\ \vfppgrp\ such that the centralizer of
every torsion element in $G$ is finite. Then $G$ is a \fpp\
product of subgroups which are finite or \fpp. \et

This is a rather surprising result from a group theoretic point of
view, since the theorem does not hold for abstract groups (as well
as for profinite groups): an easy counter example is given in
Section 5. However, from a Galois theoretic point of view it is not
so surprising. Indeed, the finite centralizer condition for
torsion elements arises naturally in the study of absolute Galois
groups. In particular, D.\,Haran \pcite{H 93} (see also I.\,Efrat
in \pcite{E 96} for a different proof) proved the above theorem
for the case when $G$ is an extension of a free pro-$2$ group with
a group of order $2$.

\medskip
The proof of Theorem 1 explores a connection between $p$-adic
representations of finite $p$-groups and virtually free pro-$p$
groups, which gives a new approach to study  virtually free
pro-$p$ groups. This connection enables us to use the following
beautiful result:

\bt{thm2-weiss}[(\pcite{weiss} A.\,Weiss)] Let $G$ be a finite $p$-group, $N$ a
normal subgroup of $G$ and let $M$ be a finitely generated
$\Z_p[G]$-module. Suppose that $M$ is a free $N$-module and that $M^N$ is a
permutation lattice for $G/N$. Then $M$ is a permutation lattice
for $G$. \et

Here $M^N$ means the fixed submodule for $N$, and a {\em
permutation lattice} for $G$ means a direct sum of $G$-modules,
each of the form $\Z_p[G/H]$ for some subgroup $H$ of $G$.

The connection to  representation theory cannot be used in a
straightforward way, however. Indeed, if one factors out the
commutator subgroup of a free open normal subgroup $F$ then the
obtained $G/F$-module would, in general, not satisfy the hypothesis
of Weiss' theorem. In order to make representation theory work, we use
 \pp\ \HNNext s to embed $G$ into a rather special \vfppgrp\ $\tilde G$,
in which, after factoring out the
commutator of a free open normal subgroup, the hypotheses of Weiss'
theorem are satisfied. With its aid we prove
Theorem 1 for $\tilde G$ and apply the Kurosh subgroup theorem to
deduce the result for $G$.

We use notation for profinite and pro-$p$ groups from \pcite{RZ1 2000}.

\section{Preliminary results}

\bt{Heller-Reiner} ((2.6) Theorem in \pcite{heller-reiner}) Let
$G$ be a group of order $p$ and $M$ a $\Z_p[G]$-module, free as a
$\Z_p$-module. Then
$$M=M_1\oplus M_p\oplus M_{p-1}$$ such that $M_p$ is a free
$G$-module, $M_1$ is a trivial $G$-module and on $M_{p-1}$ the
equality $1+c+\cdots+c^{p-1}=0$ holds, where $G=\langle
c\rangle$.\et

Let $G$ be a $p$-group. A {\em permutation lattice} for $G$ means
a direct sum of $G$-modules, each of the form $\Z_p[G/H]$ for some
subgroup $H$ of $G$. A permutation lattice will be also called
{\em $G$-permutational} module.

\medskip
If $G$ is of order $p$ then \preftheorem{Heller-Reiner} implies
that $M$ is permutational lattice if and only if $M_{p-1}$ is missing in the
decomposition for $M$ if and only if $M/(g-1)M$ is torsion free for $1\neq
g\in G$.

\bl{gp-rg} For any finite group $N$, integral domain $R$ and \fg\
free $R[N]$-module $M$ the map $\phi:M\to M^N$ defined by
$\phi(m):=\sum_{n\in N}nm$ is an \epi\ with kernel $JM$, where $J$
is the augmentation ideal in $R[N]$. \el

\begin{proof} The map is well-defined, since $\sum_{n\in N}n$ belongs to
the centre of $R[N]$ and for the same reason $JM$ is contained in
the kernel of $\phi$. Present $M=R[N]\otimes_{R[N]}L$ with $L$ a
free $R$-module, then when $m=\sum_{n\in N}n\otimes l(n)$ with
$l(n)\in L$ belongs to $M^N$, it means that all $l(n)$ are equal,
so that  $m=\phi(1_N\otimes l(n))$. Hence $\phi$ is an \epi. If
$m=\sum_{x\in N}x\otimes l(x)\in\ker\of\phi$, then mod $JM$ it is
of the form $\sum_{x\in N}1_N\otimes l(x)$ and therefore
$\sum_{x\in N}l(x)=0$ must hold, i.e., $m\in JM$. \end{proof}

\brem{$M_N$} When applying \preftheorem{thm2-weiss}, in light of
\preflemma{gp-rg}, we usually shall check the hypothesis for $M_N$
instead of $M^N$.\erem

\medskip
We shall need the following connection between free decompositions and
$\Z_p$-representations for \fpp\ by $C_p$ groups.

\bl{modules-gen} Let $G$ be a split extension of a free pro-$p$
group $F$ of finite rank by a group of order $p$. Then
\begin{rmenumerate}
\item (\pcite{Sch 99}) $G$ has a free decomposition
 $G=\left(\coprod_{a\in A}C_a\times H_a\right)\amalg H$, with $C_a\cong C_p$ and $H_a$ and $H$ \fpp.
\item Set $M:=F/[F,F]$. Fix $a_0\in A$ and a generator $c$ of
$C_{a_0}$. Then conjugation by $c$ induces an action of $C_{a_0}$
upon $M$. The latter module decomposes in the form
$$M=M_1\oplus M_p\oplus M_{p-1}$$
such that $M_p$ is a free $\gp c$-module, on $M_{p-1}$ the
equality $1+c+\cdots+c^{p-1}=0$ holds, and $c$ acts trivially on
$M_1$.

Moreover, the ranks of the three $G/F$-modules satisfy
  $\rank(M_p)=\rank(H)$, $\rank(M_{p-1})=|A|-1$, and, $\rank(M_1)=\sum_{a\in A}\rank(H_a)$.

  In particular, $M$ is $G/F$-permutational if and only if $|A|=1$.
\end{rmenumerate}
\el

\begin{proof} (i) is 1.1 Theorem in \pcite{Sch 99}.
For proving (ii), first pick $a_0\in A$, second, for each $a\in A$
a generator $c_a$ of $C_a$, and put $c_{a_0}:=c$.
We claim that
$$F=\left(\coprod_{a\in A}H_a\right)\amalg \left(\coprod_{j=0}^{p-1}H^{c^j}\right)
\amalg\left(\coprod_{a\in A\setminus\{a_0\}}\coprod_{j=0}^{p-1}\gp{c_ac\inv}^{c^j}\right).$$
Indeed, consider the \epi\ $\phi:G\to C_p$ with $H$ and all $H_a$ in the kernel and sending
each generator $c_a$ of $C_a$ to the generator of $C_p$. Then clearly $F=\ker\phi$ equals
$\gp{H^{c^j},H_a,(c_ac\inv)^{c^j}\mid a\in A,j=0,\ldots,p-1}$.
This shows that
$$\rank\of F\le
\sum_{a\in A}\rank{H_a}
+p\,\rank\of H
+(p-1)(|A|-1).
$$
On the other hand, one can use the \pp-version of the \KST, Theorem 9.1.9 in
\pcite{RZ1 2000} applied to $F$ as an open subgroup of $G$, to see that
$$F=\left(\coprod_{a\in A}H_a\right)\amalg \left(\coprod_{j=0}^{p-1}H^{c^j}\right)\amalg U$$
with $U$ a \fpp\ subgroup of $F$ having
$\rank\of U=1+|A|p-|F\backslash G/H|-\sum_{a\in A}|F\backslash G/
(H_a\times C_a)|=1+|A|p-p-|A|=(|A|-1)(p-1)$. It shows the validity of the claimed free decomposition of
$F$.

Factoring out $[F,F]$ yields the desired decomposition -- the
images of the three free factors. Finally, $M_{p-1}$ appears as
follows: writing $f_a:=c_ac\inv$ a straight forward calculation
yields  the equality $f_a^{c^{p-1}}f_a^{c^{p-2}}\cdots f_a^cf_a=1$
for every $a\in A$, which, in additive notation, reads 
$(c^{p-1}+c^{p-2}+\cdots+c+1)\bar f_a=0$.

\end{proof}

\bc{bases} If for each $a\in A$ a basis $B_{a}$  of $H_a$ is given
and $B$ is any basis of $H$, then $\bigcup_{a\in A}B_a[F,F]/[F,F]$
is a basis of $M_1$ and $B[F,F]/[F,F]$ a basis of the $G/F$-module
$M_p$. A basis of $M_{p-1}$ is given by $\{c_ac_{a_0}\inv\mid a\in
A, a\neq a_0\}[F,F]/[F,F]$.\ec

\bl{fg-cfg}
Every \fg\ \vfppgrp\ has, up to conjugation, only a finite number
of finite subgroups.
\el

\begin{proof} Suppose that the lemma is false and that $G$ is a \cex\ possessing a
normal \fpp\ subgroup $F$ of minimal possible index. When
$H$ is a maximal open subgroup  of $G$
with $F\le H$ then, as $\indx HF<\indx GF$, the proper subgroup
$H$ satisfies the conclusion of the lemma, and so there are, up to
conjugation, only finitely many finite subgroups of $G$, contained in $H$.
Hence, in order to be a \cex, $G$ must be of the form $G=F\rtimes K$ for a
finite subgroup $K$ of $G$ and, as $G$ contains only finitely many such subgroups $H$,
the proof is finished, if we can show
that up to conjugation, there are only finitely many finite
subgroups $L\cong K$ in $G$. Let $t$ be a central element of order
$p$ in $K$ and consider $G_1:=F\rtimes \langle t\rangle$.
Certainly $G_1$ is \fg. Hence, as a consequence of
\preflemma{modules-gen}(i), $G_1$ satisfies the conclusion of the
lemma, and so, $G>G_1$. Next observe that any finite subgroup
$L\cong K$ of $G$ containing some torsion element $s\in G_1$ is
contained in $C_G(s)$. By 1.2 Theorem in \pcite{Sch 99}, $C_F(s)$
is a free factor of $F$ and therefore, since $F$ is \fg, $C_F(s)$
is \fg\ as well, and so is $C_G(s)$. Let bar denote passing to the
quotient mod the normal subgroup $s$ of $C_G(s)$. Then
$\indx{\ol{C_G(s)}}{\ol{C_F(s)}}<\indx GF$, so that $\ol{C_G(s)}$
contains only finitely many conjugacy classes of
maximal finite subgroups. Since the
centralizers of conjugate elements are conjugate, $G$ can,
up to conjugation,
contain only finitely many maximal finite subgroups,
a \cdc. \end{proof}

\bl{maxsgrps} Let $G$ be \vfpp\ and $C_F(t)=\ugp$ for every
torsion element $t\in G$. Then any pair of distinct maximal finite
subgroups $A,B$ of $G$ has trivial intersection. \el

\begin{proof} Suppose that the lemma were false. Then one can pick maximal
finite subgroups $A$ and $B\neq A$ such that $1\neq C := A\cap B$
is of  maximal possible cardinality. Then $C$ is a finite normal
subgroup of $L:=\gp{N_A(C),N_B(C)}$, so the latter is itself
finite, since $N_G(C)$ must be finite (by Lemma 9.2.8 in
\pcite{RZ1 2000} a finite normal subgroup of a pro-$p$ group
intersects the centre non-trivially). On the other hand, one must
have $L\cap A=C$ due to the maximality assumption on the
cardinality of pairwise intersections of maximal finite subgroups.
Since $C<A$ one arrives at the \cdc\ $C<N_A(C)\le L\cap C=C$.
\end{proof}

We shall frequently use also the following results about virtually
\fppgrp s and \fpp\ products.

\bp{crelle}[(\pcite{Z 2003}, Proposition 1.7)] Let $G$ be a
virtually free pro-$p$ group. Then

\begin{rmenumerate}
\item $G/\torgp{G}$ is free pro-$p$;

\item $\tor G$ maps onto $\tor{G/\torgp{G}}$ under the canonical
epimorphism  $G\longrightarrow G/\torgp{G}$.
\end{rmenumerate}
\ep

\bt{conjugacy}[(\pcite{RZ1 2000}, Theorems 9.1.12 and 9.5.1)]
Let $G=\coprod_{i=1}^n G_i$ be a free
profinite (pro-$p$) product. Then $G_i\cap G_j^g=1$ for either
$i\neq j$ or $g\not\in G_j$.

Every finite subgroup of $G$ is conjugate to a subgroup of
a free factor.
\et

\section{HNN-embedding}
We introduce a notion of a pro-$p$ \HNNgrp\ as a
generalization of pro-$p$ \HNNext\ in the sense of
\pcite{RZ2}, page 97. It also can be defined as a sequence of pro-$p$
\HNNext s. During the definition to follow, $i$ belongs to a
finite set $I$ of indices.

\bd{HNN-grp} Let $G$ be a pro-$p$ group and $A_i, B_i$ be
subgroups of $G$ with isomorphisms $\phi_i:A_i\longrightarrow
B_i$. The pro-$p$ \HNNgrp\  is then a pro-$p$ group
$\HNN(G,A_i,\phi_i,z_i)$ having presentation
$$\HNN(G,A_i,\phi_i,z_i)=\langle G, z_i\mid rel(G), \forall a_i\in A_i:\ \
a_i^{z_i}=\phi_i(a_i)\rangle.$$ The group $G$ is called the {\em
base group}, $A_i,B_i$ are called {\em associated subgroups} and
$z_i$ are called the {\em stable letters}. \ed

For the rest of this section let $G$ be a \fg\ \vfppgrp, and
fix an open \fpp\ normal subgroup $F$ of $G$ of minimal index.
Also suppose that $C_F(t)=\ugp$ for every torsion element $t\in
G$. Let $K:=G/F$ and form $G_0:=G\amalg K$. Let
$\psi:G\to K$ denote the canonical projection.
It extends to an \epi\ $\psi_0:G_0\to K$, by sending $g\in G$ to $gF/F\in K$ and each $k\in K$ identically to $k$, and using the universal
property of the \fpp\ product. Remark that the kernel
of $\psi_0$, say $L$, 
is an open subgroup of $G_0$ and, as $L\cap G=F$ and $L\cap K=\ugp$,
as a consequence of the \pp\ version of the \KST, Theorem 9.1.9 in
\pcite{RZ1 2000}, $L$ is \fpp.
Let $I$ be the set
of all $G$-conjugacy classes of maximal finite subgroups of $G$.
Fix, for every $i\in I$, a finite subgroup $K_i$ of
$G$ in the $G$-conjugacy class $i$. We define a pro-$p$ \HNNgrp\
by considering first  $\tG_0:=G_0\amalg F(z_i\mid i\in I)$ with $z_i$
constituting a free set of generators, and then taking the normal
subgroup $R$ in $\tG_0$ generated by all elements of the form
$k_i^{z_i}\psi(k_i)\inv$, with $k_i\in K_i$ and $i\in I$. Finally
set
$$\tG:=\tG_0/R,$$
and note that it is an \HNNgrp\ $\HNN(G_0, K_i,\phi_i,z_i)$,
where $\phi_i:=\psi_{|K_i}$,  $G_0$ is the base group, the $K_i$
are associated subgroups, and the $z_i$ form a set of stable
letters in the sense of \prefdef{HNN-grp}.

\medskip
Let us show that $\tG$ is \vfpp. The above \epi\
$\psi_0:G_0\longrightarrow K$ extends to $\tG\longrightarrow K$
by the universal property of the HNN-extension, so $\tG$ is a
semidirect product $\tilde F\rtimes K$ of its kernel $\tilde F$
with $K$. By Lemma 10 in \pcite{HZ 07}, every open torsion free
subgroup of $\tG$ must be \fpp, so $\tilde F$ is free pro-$p$.

\medskip

The objective of the section is to show that the centralizers of
torsion elements in $\tilde G$ are finite.

\bl{centr} Let $\tilde G= \HNN(G_0,K_i,\phi_i,z_i)$ and $\tilde F$
be as explained.   Then $C_{\tilde F}(t)=1$ for every torsion
element $t\in \tilde G$. \el

\begin{proof} There is a standard pro-$p$ tree $S:=S(\tilde G)$
associated to $\tilde G:= \HNN(G_0,K_i,\phi_i,z_i)$ on which
$\tilde G$ acts naturally such that the vertex stabilizers are
conjugates of $G_0$ and each edge stabilizer is a conjugate of
some $K_i$ (cf. \pcite{RZ2} and \S 3 in \pcite{Z-M 90}).

\medskip

\noindent{\em Claim:
Let $e_1,e_2$ be two edges of $S$ with a common vertex
$v$. Then the intersection of the stabilizers $\tilde
G_{e_1}\cap\tilde G_{e_2}$ is trivial.
}

\medskip

\claimproof By translating $e_1,e_2, v$ if necessary we may assume
that $G_0$ is the stabilizer of $v$. Then, up to orientation, we
have two cases:

1) $v$ is initial vertex of $e_1$ and $e_2$. Then $\tilde
G_{e_1}=K_i^g$ and $\tilde G_{e_2}=K_j^{g'}$ with $g,g'\in G_0$
and either $i\neq j$ or $g\not\in K_ig'$. Suppose that $K_i^g\cap
K_j^{g'}\neq\ugp$. Then, since $G_0=G\amalg K$, we may apply
\preftheorem{conjugacy}, in order to deduce the existence of
$g_0\in G_0$ with $K_i^{gg_0}\cap K_j^{g'g_0}\le G$. Now apply
\preflemma{maxsgrps}, in order to deduce the \cdc\ $i=j$ and
$gg_0\in K_ig'g_0$. So we have $K_i^g\cap K_j^{g'}=\ugp$, as
needed.
\medskip

2) $v$ is the terminal vertex of $e_1$ and the initial vertex of
$e_2$. Then $\tilde G_{e_1}=K^g$ and $\tilde G_{e_2}=K_i^{g'}$ for
$g,g'\in G_0$ so they intersect trivially by the definition of
$G_0$ and \preftheorem{conjugacy}. So the Claim holds.

\bigskip

Now pick a torsion element $t\in \tilde G$ and $f\in\tilde F$ with
$t^f=t$. Let $e\in E(S)$ be an edge stabilized by $t$. Then $fe$
is also stabilized by $t$ and, as by Theorem 3.7 in \pcite{RZ2},
the fixed set $S^t$ is a subtree, the path $[e,fe]$ is fixed by
$t$ as well. By the above then $fe=e$ contradicting the freeness
of the action of $\tilde F$ on $E(S)$.

\end{proof}
\section{Proof of the main result}

\bp{model1} Let $G$ be a semidirect product of a \fppgrp\ $F$ of
finite rank with a $p$-group\ $K$ such that every finite subgroup is
conjugate to a subgroup of $K$. Suppose that $C_F(t)=\ugp$ holds for
every torsion element $t\in G$. Then $G=K\amalg F_0$ for a \fpp\
factor $F_0$. \ep

\begin{proof}
Suppose that the proposition is false. Then there is a \cex\ with $K$ having
minimal order. When $K\cong C_p$, then by
\preflemma{modules-gen}(i) $G=\left(\coprod_{i\in I}C_i\right)\amalg H$ with
$I$ a finite set, all $C_i$ of order $p$ and $H$ \fpp.
By the assumptions and \preftheorem{conjugacy} there is a single conjugacy
class of finite subgroups, i.e., $|I|=1$, so that $G$ would not be a \cex.
Therefore  $K$ is of order $\ge
p^2$.

Let $H$ be any maximal subgroup of $K$. Then $F\rtimes H$
satisfies the premises of the proposition and hence $F\rtimes H$
is of the form $H\amalg F_1$ for some free factor $F_1$. Let us
denote by bar passing to the quotient mod $(H)_G$. As
$(H)_G=\torgp{FH}$ by \prefprop{crelle}(i) $\bar F$ is \fpp.
\preflemma{modules-gen}(i) shows that $\bar G\cong \coprod_{i\in
I}\left(C_i\times C_{\bar F}(C_i)\right)\amalg F_0$ with $I$
finite and $F_0$ a free factor of $\bar F$. Now by
\prefprop{crelle}(ii)  $\overline{\tor{\bar G}}=\tor{\bar G}$,
and therefore, every torsion element in $\bar G$ can be lifted to
a conjugate of an element in $K$. Hence $I$ consists of a single
element, so that \be{FH} \bar G=(\bar K\times C_{\bar F}(\bar
K))\amalg F_0. \ee In the sequel we shall use
\preflemma{modules-gen}(ii) a couple of times.
 Consider $M:=F/F'$ as
a $K$-module and let $J$ denote the augmentation ideal of
$\Z_p[H]$. Since $F\rtimes H=H\amalg F_1=\left(\coprod_{h\in
H}F_1^h\right)\rtimes H$, $H$ acts by permuting the free factors
$F_1^h$, so that $M$ is a free $H$-module.  Passing in \prefeq{FH}
to the quotient mod the commutator subgroup of $\bar F=(C_{\bar
F}(\bar K),F_0)_{\bar G}$, using \preflemma{modules-gen}, one can
see that $M/JM$ is  a $\bar K$-permutation lattice. Then  an
application of \preftheorem{thm2-weiss} together with
\prefrem{$M_N$} shows that $M$ itself is a $K$-permutation
lattice.

We shall  show that $M$ is a free $K$-module. Indeed, if any of
the summands is not free, a proper subgroup of $K$, say $S$, acts
trivially there. Since $M$ is a free $H$-module, conclude that
$S\cap H=\ugp$, and from this that $S$ is of order $p$. Let us
show that $G_1:=F\rtimes S$ satisfies the premises of the
proposition. Certainly $C_F(t)=\ugp$ for every torsion element
$t\in G_1$. Pick $x\in \tor{G_1}$. There is $k\in K$ and $f\in F$
with $x=k^f$. Since $k\in (FS)\cap K$ deduce $k\in S$. So there is
a single conjugacy class of finite subgroups in $G_1$. But then,
considering the natural homomorphism from $F\rtimes S$ to
$M\rtimes S$ and, observing the \minass\ on $|K|>|S|=p$, so that
$F\rtimes S=S\amalg F_S$ for some \fppgrp\ $F_S$,
 one finds as an application of \preflemma{modules-gen} that the
decomposition of $M$ cannot have direct summands, on which $S$
acts trivially, a \cdc. Since $M$ is a $K$-permutational lattice,
it is $S$-permutational as well and so  cannot contain $p-1$
blocks, so that $M$ is a free $K$-module.

Consider $\tG:=K\amalg \tF_0$ with $\tF_0\cong \torfactor G$. By
\prefprop{crelle}  $\torfactor G$
is free pro-$p$, so we can  fix a section $F_0$ of $\torfactor G$
inside $G$, and define a \homo\ $\phi:\tG\to G$ by first sending $K$ to
$K$ and $\tF_0$ onto $F_0$ and then, using the universal property
of the free product, extending it to $\tG$.
By assumption all torsion elements are, up to conjugation, contained in $K$, showing
that $\phi$ is an \epi.
By the above the kernel of $\phi$ must be
contained in $[\tF,\tF]$. In particular, since the group is \fg, one
has $\tF\cong F$, since both groups are \fpp. Since
$K\cap\ker\phi=\ugp$, conclude that $\phi$ is an isomorphism, as
claimed. \end{proof}

\begin{proof}[of \preftheorem{main}:]
\preflemma{fg-cfg} shows that $G$ can have only a finite
number of conjugacy classes of maximal finite subgroups. Therefore
one can form $\tG$ as described before \preflemma{centr}, in order
to embed $G$ such that $\tG$ is both, \fg, and, has finite
centralizers of its finite subgroups, and, moreover, has a single
conjugacy class of maximal finite subgroups.  By \prefprop{model1}
the group $\tG$ is of the form $\tG=K\amalg F_0$ where $K$ is
finite and $F_0$ is \fpp. Since $G$ is a \fg\ \pp\ subgroup of
$\tilde G$, the \KST\ in \pcite{HR 87} implies that
$G$ must have indeed the form as claimed.
\end{proof}

\section{An example}

We give an example of a virtually free profinite group that satisfies
the centralizer condition of the main theorem but does not satisfy
its conclusion. Note that the same example is valid for abstract
groups.

\bl{cex} Let $A\cong B=S_3$ be the symmetric group on a
$3$-element set and $C:=C_2$. Form the amalgamated
free profinite product $G=A\amalg_CB$, where $C$ identifies with
given $2$-Sylow subgroups in $A$ and $B$ respectively.

Then for every torsion element $t\in G$ its centralizer is finite. However,
$G$ cannot be decomposed as a free profinite product with some factor finite.
\el

\begin{proof} It is easy to see that $G$ can be presented in the form
$G=N\rtimes C_2$, with $N\cong C_3\amalg C_3$ and $C_2=\gp\alpha$
acting by inverting the generators of the two factors. Then the
structure of $N$, in light of \preftheorem{conjugacy}, shows that
no element of order $3$ can have an infinite centralizer. For
establishing the first statement of the Lemma, it will suffice to
show that all involutions in $G$ are conjugate, and that $\alpha$
acts without fixed points upon $N=\gp{a,b}$, where $a,b$ are
generators of cyclic free factors of order $3$. As $G$ is the
fundamental group of the graph of groups $\xy\POS(0,0)
*{\bullet}\POS(-2,2) *{^{A}} \POS(20,0) *{\bullet}\POS(22,2)
*{^{B}} \POS(0,0)\ar(20,0)^{C}
\endxy $,  Theorem 5.6 on page 938 in \pcite{Z-M 90} shows that every
involution is conjugate to an involution in one of the vertex
groups $A$ or $B$. As $A$ and $B$ both have a single conjugacy
class of involutions and the latter contains $C$, the first
observation holds. Since, by Theorem 9.1.6 in \pcite{RZ1 2000},
$N'$ is freely generated by the commutators $[a^i,b^j]$ with
$i,j\in \{1,2\}$, one can see that $\alpha$ permutes them without
fixed points, so that $N'\rtimes\gp\alpha$ is isomorphic to
$F(x,y)\amalg C_2$ with $F(x,y)$ a free profinite group. Thus
$\alpha$ has no fixed points in $N'$ and, as an easy consequence,
none in $N$.

Suppose that $G=L\amalg K$ with $L$ finite.
Then, by
\preftheorem{conjugacy}, w.l.o.g. we can assume that  $A\le L$.
The just cited Theorem on page 938 in \pcite{Z-M 90}
shows that $A$ is a maximal finite subgroup of $G$, so that
$A=L$. Since the quotient mod the normal closure of $L$ in $G$ is
isomorphic to $K$ on the one hand and trivial by construction,
find $K=\ugp$, a \cdc. So $G$ has no finite free factor. \end{proof}

A list of remarks of a referee of a previous version of the paper led to immense improvement in
presenting some proofs.
\newcommand{\bibname}[1]{{\sc #1}}

\begin{tabular}l
W. Herfort\\
University of Technology Vienna, Austria\\
{w.herfort@tuwien.ac.at}
\end{tabular}\hfill\begin{tabular}l
P. A. Zalesski\\
University of Brasilia, Brazil\\
{pz@mat.unb.br}\end{tabular}
\end{document}